\documentclass[a4paper,10pt,twocolumn]{article}

\usepackage[english]{babel}
\usepackage[utf8]{inputenc}
\usepackage[T1]{fontenc}

\usepackage[top=1.5cm, left=1.5cm, right=1.5cm, bottom=1.5cm]{geometry}

\renewenvironment{abstract}{\bf\small {\em\ Abstract---}}{}

\usepackage{amsfonts,amssymb,amsmath,amsthm}
\usepackage{subfigure}
\usepackage{graphicx}
\usepackage[footnotesize]{caption}

\usepackage{amsfonts,amssymb,amsmath,amsthm}
\usepackage{url}
\usepackage{bm}
\usepackage{grffile}
\usepackage{siunitx}
\title{Mathematical osmosis imaging for multi-modal and multi-spectral applications in Cultural Heritage conservation}

\author{Simone Parisotto$^1$, Luca Calatroni$^2$ and Claudia Daffara$^{3}$.\\ \footnotesize 
  $^1$CCA, University of Cambridge (UK).\ 
  $^2$CMAP, \'Ecole Polytechnique  (France) .\
  $^3$University of Verona (Italy).  } 
  \date{\empty} 

\begin{document}

\maketitle

\begin{abstract} 
In this work
we present a dual-mode mid-infrared workflow \cite{Daffara2018}, for detecting sub-superficial mural damages in frescoes artworks.
Due to the large nature of frescoes, multiple thermal images are recorded. Thus, the experimental setup may introduce measurements errors, seen as inter-frame changes in the image contrast, after mosaicking.
An approach to lowering errors is to post-process the mosaic \cite{Parisotto2018} via \emph{osmosis} partial differential equation (PDE) \cite{Vogel2013,Weickert2013}, which preserves details, mass and balance the lights: efficient numerical study for osmosis on large images is proposed \cite{Calatroni2017,Parisotto2018general}, based on operator splitting \cite{Hundsdorfer}.
Our range of Cultural Heritage applications include the detection of sub-superficial voids in \emph{Monocromo} (L.\ Da Vinci, Castello Sforzesco, Milan) \cite{Daffara2015}, the light-balance for multi-spectral imaging and the data integration on the Archimedes Palimpsest \cite{Parisotto2018}. 
\end{abstract}

\section{Introduction}
\label{sec: introduction}
Nowadays, image processing is an active research field in Cultural Heritage (CH) conservation. 
In particular, multi-spectral image analysis can reveal hidden features of artworks and help specialists during the restoration, e.g.\ \cite{Asperen,Daffara11,Gavrilov,Maldague}. 
These non-destructive methods offer high-resolution images and large datasets to be processed via efficient algorithms, eventually.

In this work, we describe a new multi-modal mid-infrared (MWIR, \SIrange{3}{5}{\micro\meter}) imaging technique \cite{Daffara2018,Daffara2015} to detect sub-superficial defects in fresco walls, by combining the emitted and reflected information recorded by a thermal camera with suitable lighting sources. 
Here, multiple images are acquired and mosaicked together so as to form a large MWIR mosaic. 
However, the experimental setup may introduce measurement errors, seen in our applications as inter-frame changes in contrast (i.e.\ constant \emph{shadows}). 
Thus, errors are lowered via a recently introduced drift-diffusion equation \cite{Vogel2013,Weickert2013}, called \emph{osmosis}, already able to remove constant shadows in images.
We implemented osmosis efficiently for large images (up to 30 MP), with standard numerical splitting approaches \cite{Hundsdorfer,Calatroni2017}.

\paragraph*{Applications.}
We discuss the dual-mode workflow for the detection of sub-superficial voids and cement patches in the restoration of \emph{Monocromo}, a fresco wall by Leonardo Da Vinci (Castello Sforzesco, Milan). Results are validated by specialist restorers from \emph{Opificio delle Pietre Dure} (Florence). 

We applied the efficient osmosis model for different CH images acquired in different electromagnetic bands, where a light balance post-processing enhances the information: MWIR thermal quasi-reflectography, UV fluorescence and IR falsecolor imaging.
We also discuss the use of osmosis for a multi-spectral data-fusion on the Archimedes Palimpsest.

\paragraph{Organisation.}
In Section \ref{sec: multimodal}, we describe the core idea of the multi-modal MWIR workflow for the detection of mural defects.
In Section \ref{sec: osmosis}, we recall the osmosis equation and the splitting approach for an efficient numerical solution of the light balancing problem.
In Section \ref{sec: otherapp}, we present other applications of the osmosis model in different CH imaging problems.

\section{Dual-mode mid-infrared imaging}
\label{sec: multimodal}
In \cite{Daffara2012}, a new MWIR methodology called \emph{Thermal Quasi-Reflectography} (TQR) is proposed for revealing more sub-superficial mural features than traditional approaches.
At the core of TQR there is the observation that an object, at room temperature, emits in the MWIR spectrum only 1\% of its thermal energy: thus, its (quasi-)reflected radiation $\rho$ can be recorded by a MWIR thermal camera and suitable light sources.

In \cite{Daffara2018}, this idea is evolved into a dual-mode (reflected plus emissivity) MWIR infrared approach. 
Here, TQR is the first step of the workflow since it provides the emissivity $\varepsilon=1-\rho$ via the simplified Kirchhoff's law. 
Such  $\varepsilon$ is then plugged in the recording of a thermal cooling down sequence in the emissivity mode.
Both steps are performed by the same thermal device in a fixed geometry lighting setup, tuned safely for artworks.

Thus, the dual-mode complementary information can be superimposed to unveil hidden mural features, e.g.\ different cooling down profiles. 
Also, TQR details allow to register the two aligned datasets onto a visible orthophoto, an impossible task for the emissivity-only mode due to the heat diffusion blurring.

\begin{figure}[tbhp]
\centering
\includegraphics[width=0.116\textwidth]{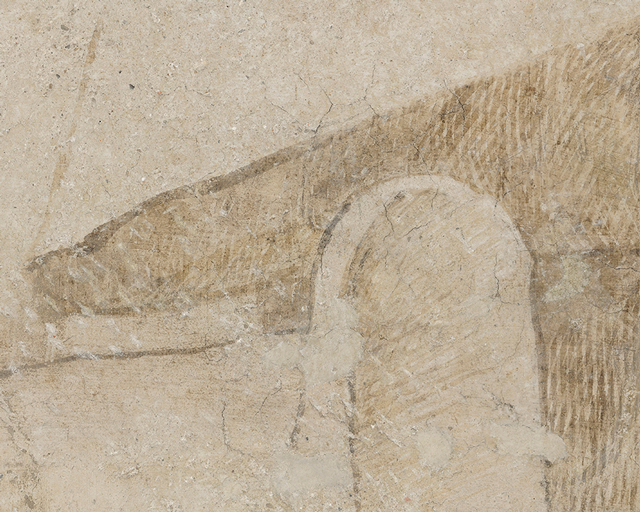}\,
\includegraphics[width=0.116\textwidth]{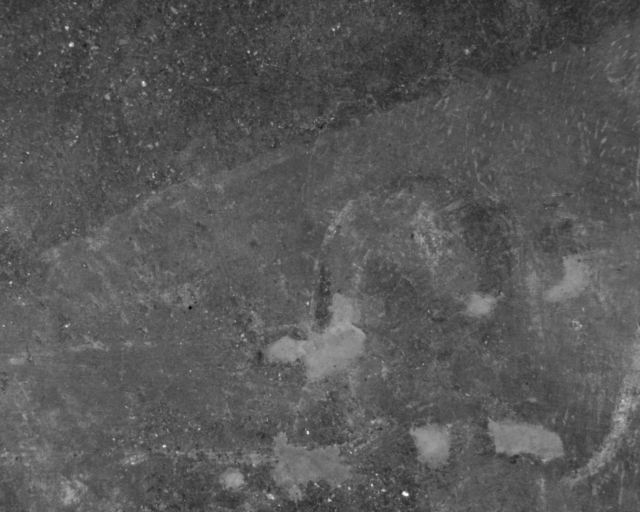}\,
\includegraphics[width=0.116\textwidth]{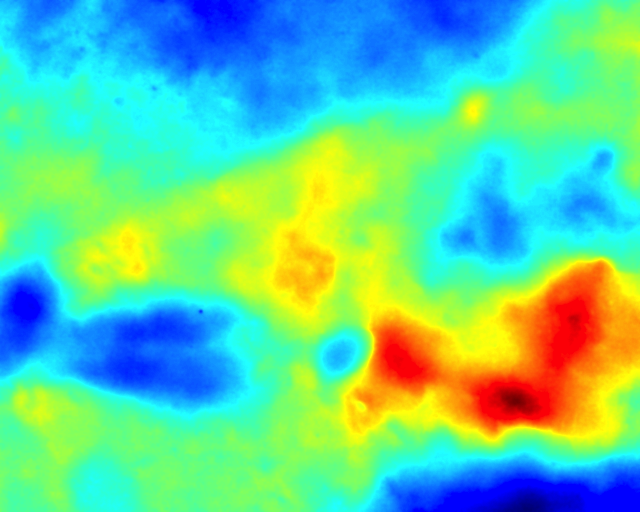}\,
\includegraphics[width=0.116\textwidth]{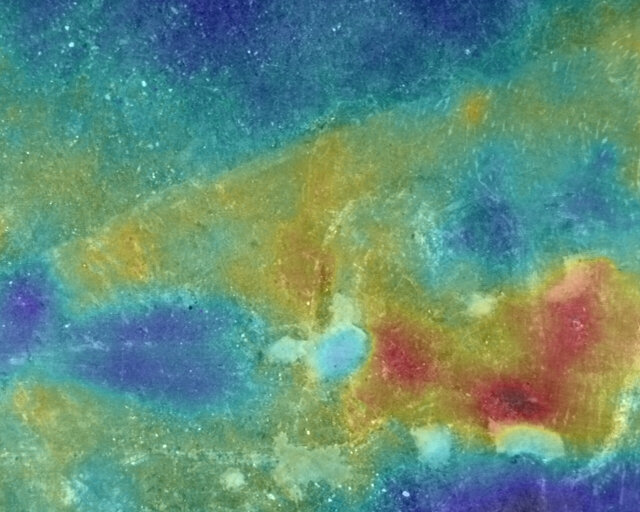}
\\
\includegraphics[width=0.116\textwidth,trim=0cm 0 0 10cm,clip=true]{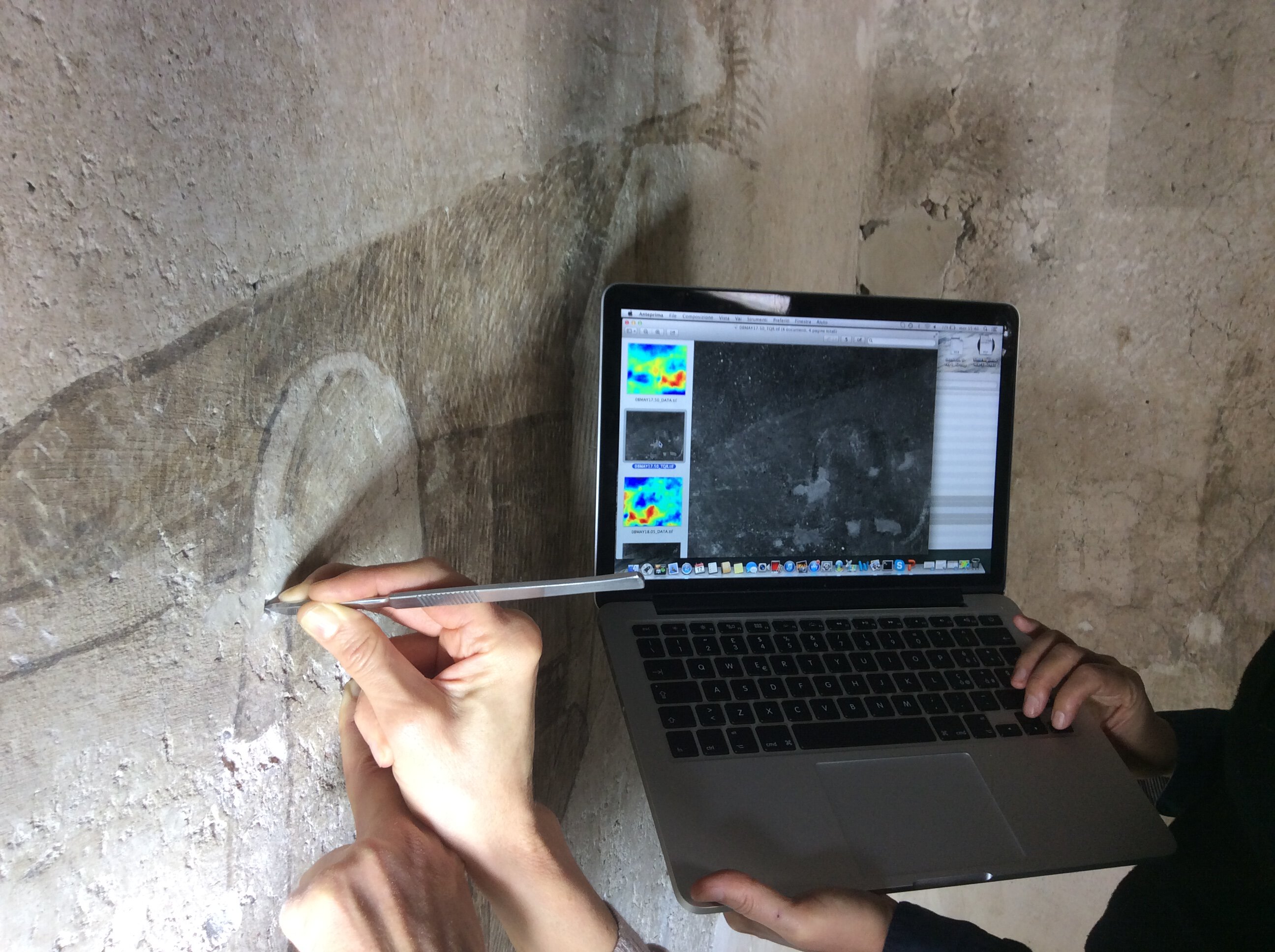}\,
\includegraphics[width=0.116\textwidth,trim=0cm 0 0cm 10cm,clip=true]{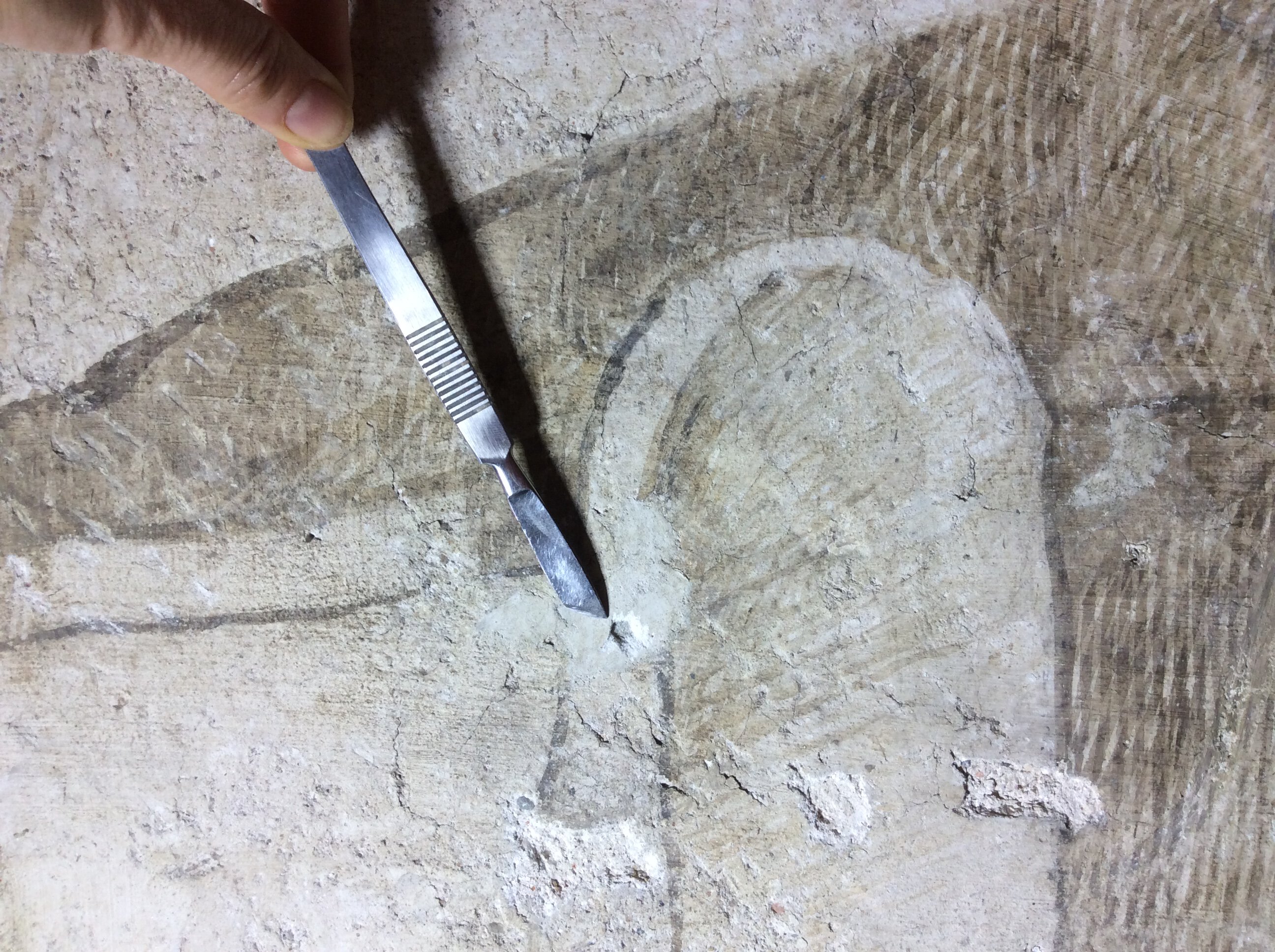}
\caption{Dual-mode MWIR imaging for discovering hidden cement patches.}
\label{fig: dual-mode}
\end{figure}

In Fig.\ \ref{fig: dual-mode} we report (top row) the visible fresco, the reflected TQR, the emitted thermal images and the fusion of the complementary data; we also show (bottom row) the discovering of sub-superficial cement patches via the dual-modal approach.

\begin{figure}[tbhp]
\centering
\includegraphics[width=0.13\textwidth]{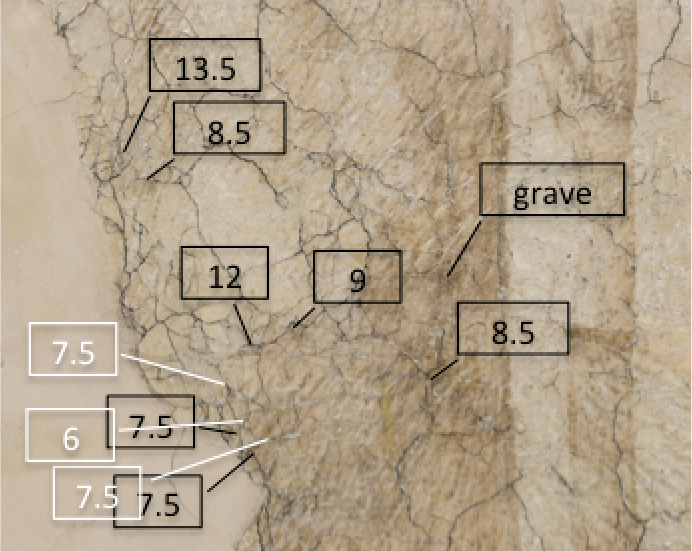}\,
\includegraphics[width=0.13\textwidth]{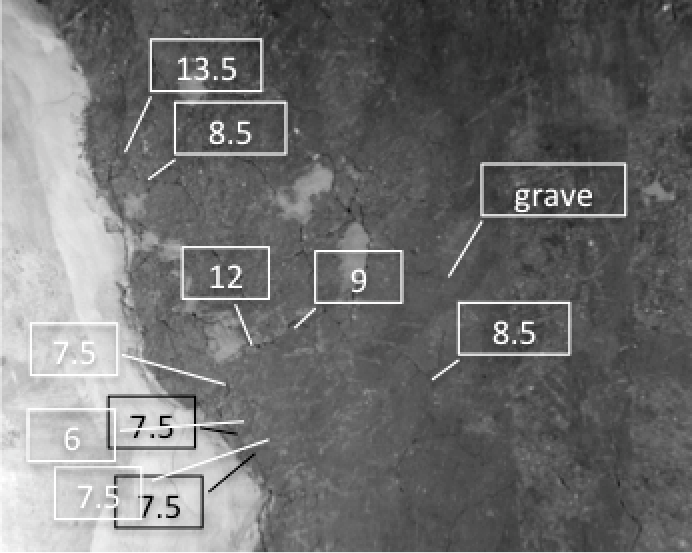}\,
\includegraphics[width=0.13\textwidth]{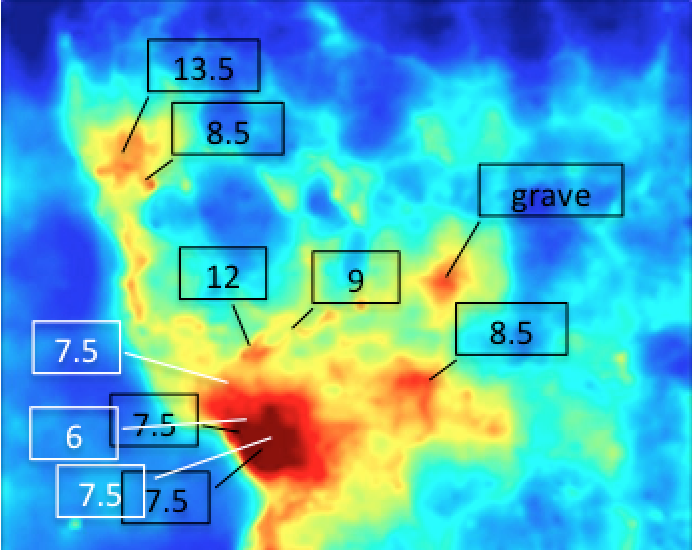}
\caption{Defect measures (in mm.) of identified damages via dual-mode setup.}
\label{fig: pin test}
\end{figure}

In Fig.\ \ref{fig: pin test}, sub-superficial voids are effectively detected: different void shapes at sub-millimetric precision are revealed and measured with pins by restorers of \emph{Opificio delle Pietre Dure}.


\section{Efficient osmosis for light balancing}
\label{sec: osmosis}
Let $\Omega\subset\mathbb{R}^2$, $u,f,v:\Omega\to\mathbb{R}$, $\mathbf{d}:\Omega\to\mathbb{R}^2$ with $\mathbf{d}=\nabla\log v$. 
The following model is called \emph{osmosis} \cite{Vogel2013,Weickert2013}:
\begin{equation}
\begin{cases}
u_t = \mathrm{div}(\nabla u - \mathbf{d} u) &\text{on } \Omega\times (0,T];\\
u(0)= f(x) &\text{on } \Omega \text{ at }t=0;\\
\langle \nabla u-\mathbf{d} u,\mathbf{n}\rangle = 0 &\text{on }\partial\Omega.
\end{cases}
\label{eq: osmosis}
\end{equation}
Interestingly, the non-symmetric drift-diffusion PDE \eqref{eq: osmosis} has nice properties, e.g.\ conservation of the mass, positivity and convergence to non-constant steady states, which makes its use appealing for imaging applications, when $\mathbf{d}$ is locally modified. 
This is the case for removing constant shadows in images, where $\mathbf{d}=0$ on the shadow edges, i.e.\ the jumps produced by sharp shadows.
Numerical solutions of \eqref{eq: osmosis} are given by solving the linear system $u_t = \mathbf{A} u$ where the 5-diagonal spatial discretisation $\mathbf{A}$ of $\mathrm{div}(\nabla u - \mathbf{d} u)$ is given by $\mathbf{A} = \mathbf{A}_{\bm{1}} + \mathbf{A}_{\bm{2}}$
with
\[\small
\begin{aligned}  
\mathbf{A}_{\bm{1}}(u) :=
&\dfrac{u_{i+1,j}-2u_{i,j}+u_{i-1,j}}{h^2} 
\\
&- \left( d_{1,i+\frac{1}{2},j}\dfrac{u_{i+1,j}+u_{i,j}}{2h}-d_{1,i-\frac{1}{2},j}\dfrac{u_{i,j}+u_{i-1,j}}{2h}\right)
\\
\mathbf{A}_{\bm{2}}(u) :=
&\dfrac{u_{i,j+1}-2u_{i,j}+u_{i,j-1}}{h^2} 
\\
&-
\left( d_{2,i,j+\frac{1}{2}}\dfrac{u_{i,j+1}+u_{i,j}}{2h}-d_{2,i,j-\frac{1}{2}}\dfrac{u_{i,j}+u_{i,j-1}}{2h} \right).
\end{aligned}
\]
From \cite[Prop.\ 2]{Vogel2013}, $\mathbf{A}$ is irreducible, with non-negative off-diagonal entries and column sums equal to 0: thus, scale-space properties are guaranteed for Implicit and Explicit (with time-step restriction) Euler \cite[Prop.\ 1]{Vogel2013}.

\begin{figure}[tbhp]
\centering
\includegraphics[width=0.13\textwidth]{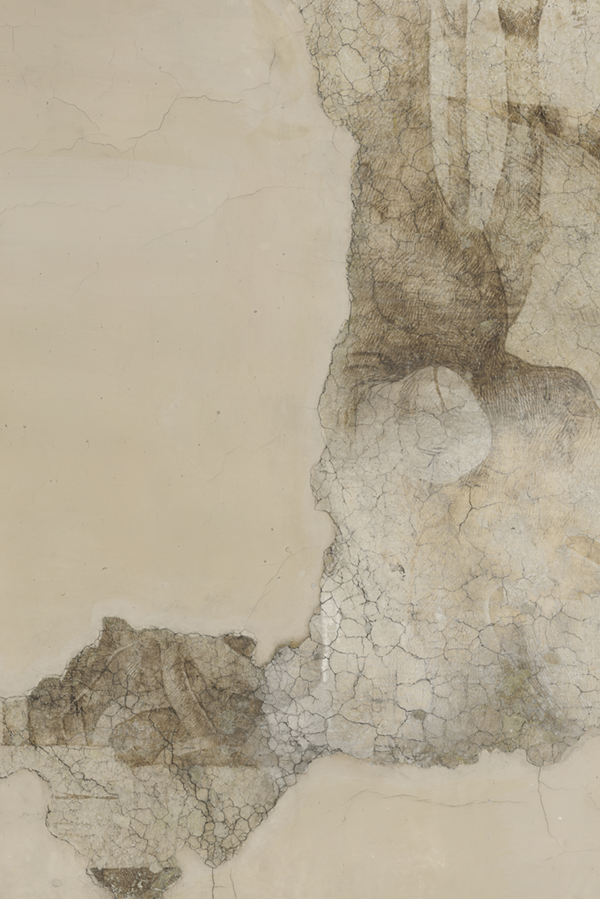}\,
\includegraphics[width=0.13\textwidth]{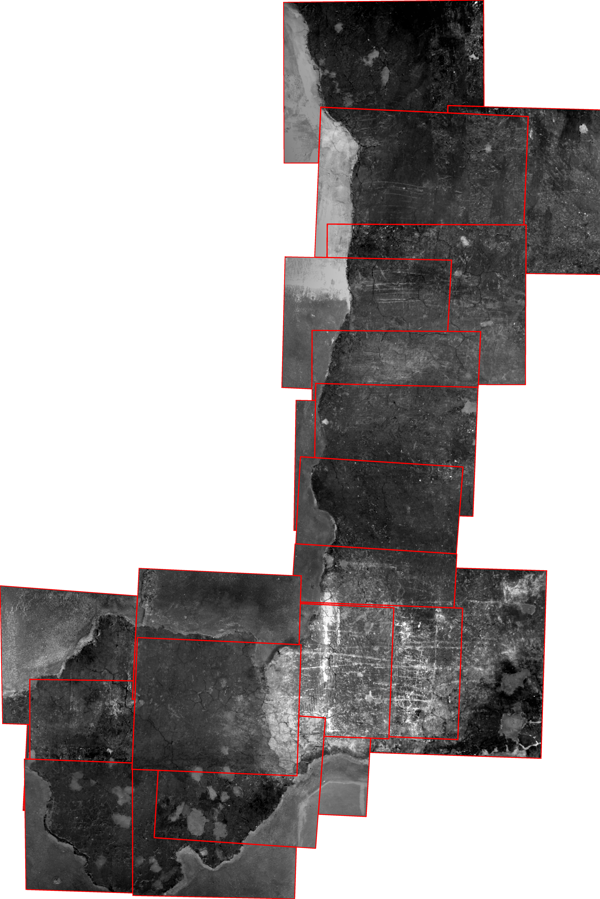}\,
\includegraphics[width=0.13\textwidth]{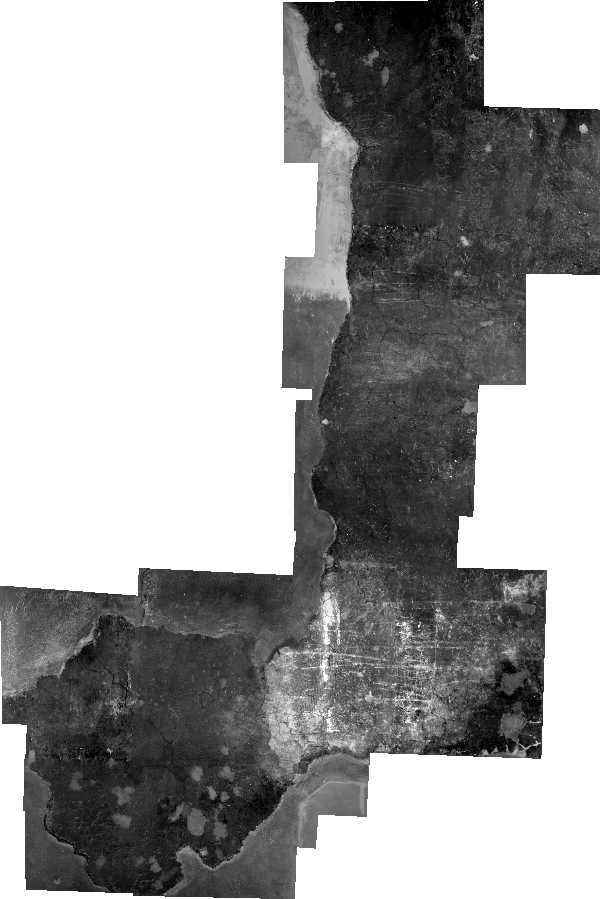}
\caption{\emph{Monocromo}, L.\ Da Vinci. Visible (left), TQR mosaicked (middle), osmosis result (right, 28 MP processed in 629 s.\ via MOS splitting)}
\label{fig: TQR balance}
\end{figure}

For large images, a way to avoid the computational bottleneck due to the large bandwidth of $\mathbf{A}$ is to consider $\mathbf{A}_{\bm{1}}$ and $\mathbf{A}_{\bm{2}}$ separately, by ADI splitting \cite{Hundsdorfer}. 
Among the others, the Multiplicative Operator Splitting (MOS) scheme preserves the scale-space properties of \eqref{eq: osmosis} and provides effective results for large time-steps $\tau>0$ \cite{Parisotto2018general}:
\begin{equation}
u^{k+1} 
= 
\prod_{\bm{n}=1}^2~\left(\mathbf{I} -\tau \mathbf{A}_{\bm{n}}\right)^{-1} u^{k}, \tag{MOS}  \label{eq:MOS}
\end{equation}
In Fig.\ \ref{fig: TQR balance} we applied the osmosis filter to a TQR mosaic obtained from the workflow described in Section \ref{sec: multimodal}. Thus, we aimed to balance the inter-frame contrasts while preserving intra-frame details: this helped restorers in the visual inspection and comparison of hidden features of the mural painting.

\section{Other CH applications}
\label{sec: otherapp}
\begin{figure}[tbhp]
\centering
\includegraphics[width=0.13\textwidth]{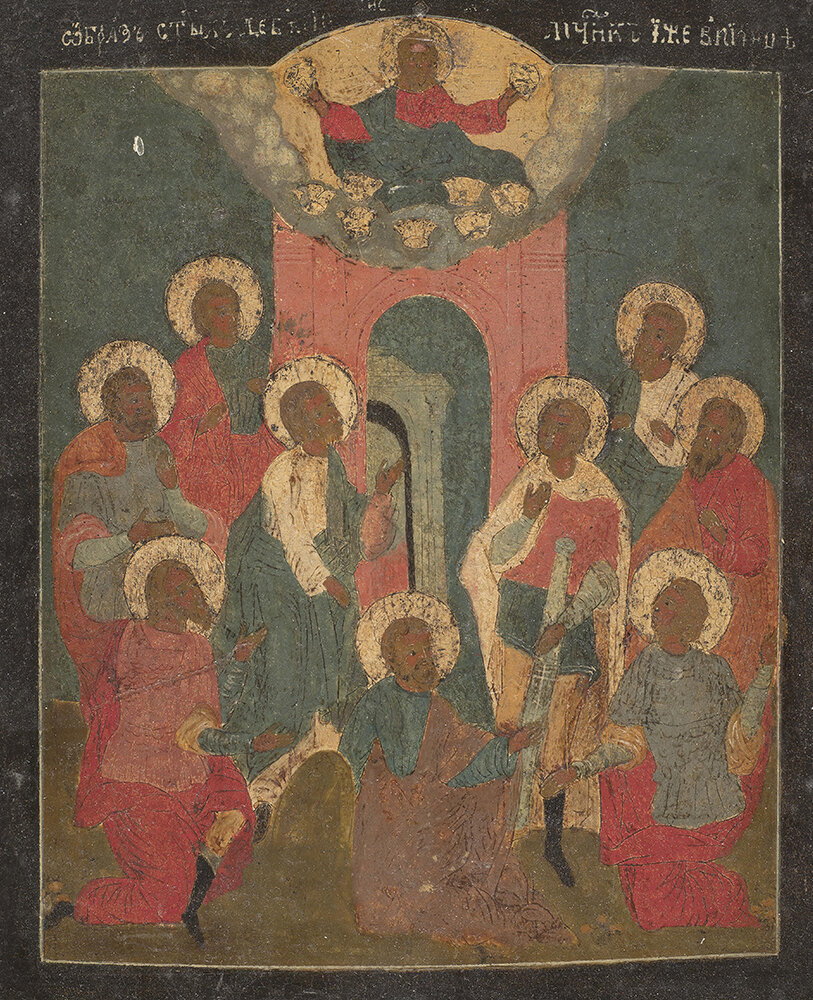}\,
\includegraphics[width=0.13\textwidth]{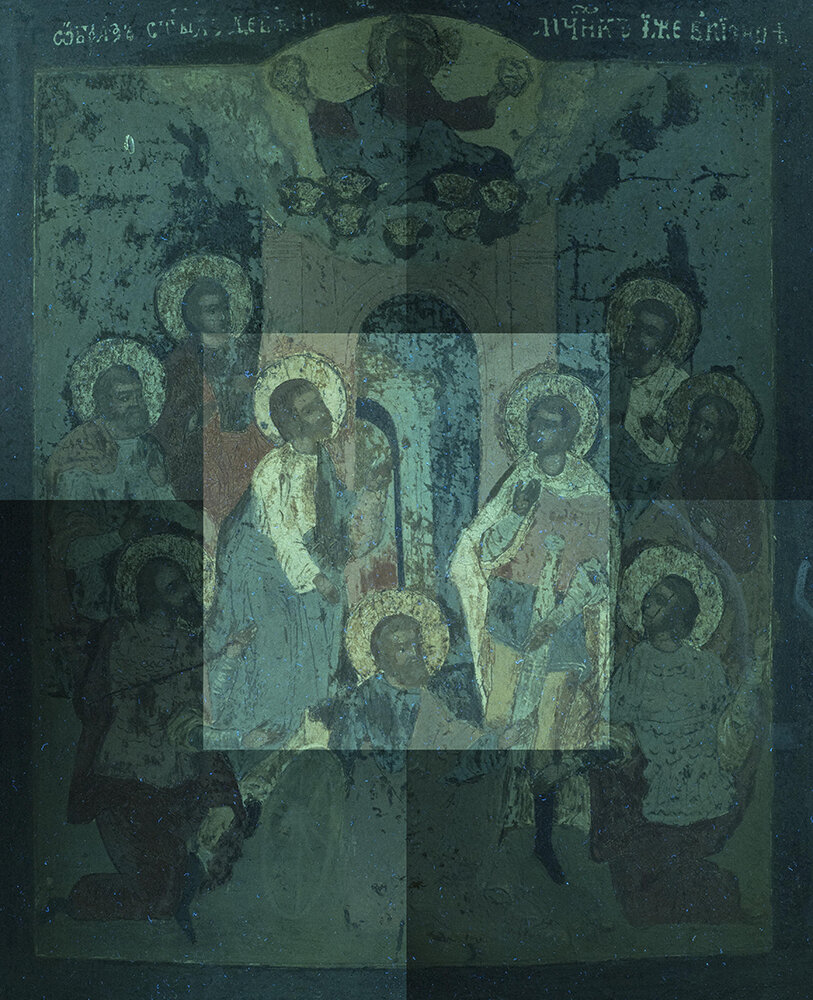}\,
\includegraphics[width=0.13\textwidth]{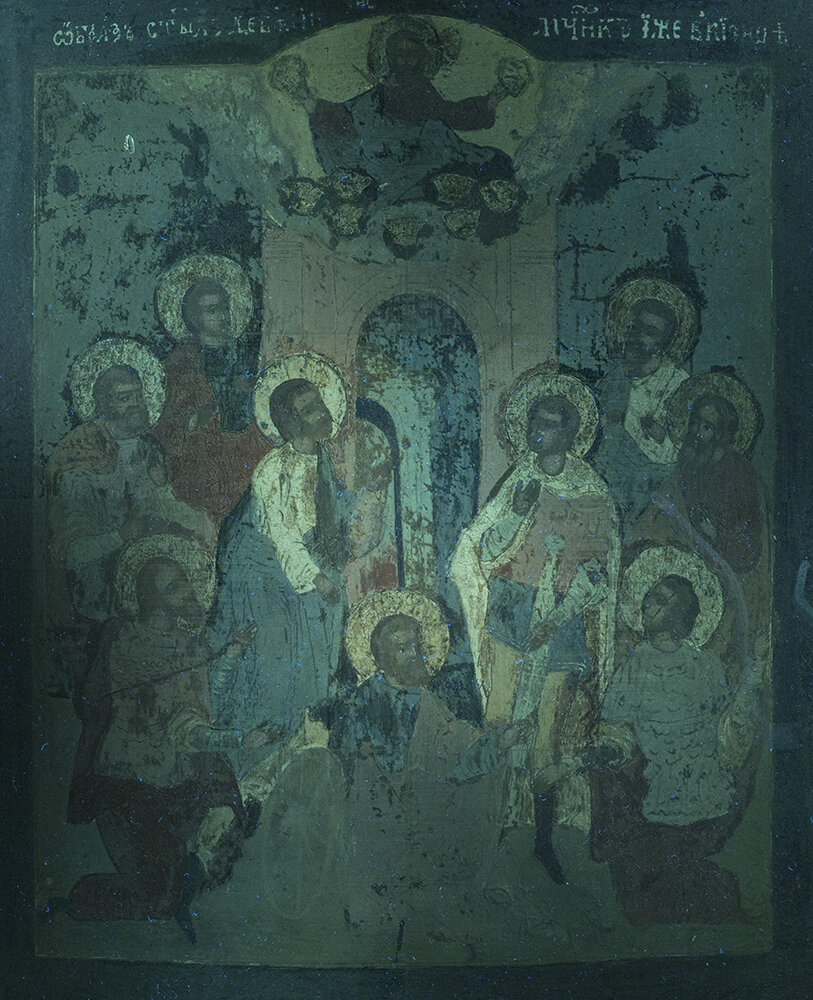}
\caption{Light balancing in UV Fluorescence. Visible (left), 5 UV Fluorescence tiles, result (right, 18 MP and 3 color channels, 1693 s. \cite{Parisotto2018})}
\label{fig: uvf}
\end{figure}

\begin{figure}[tbhp]
\centering
\includegraphics[width=0.13\textwidth]{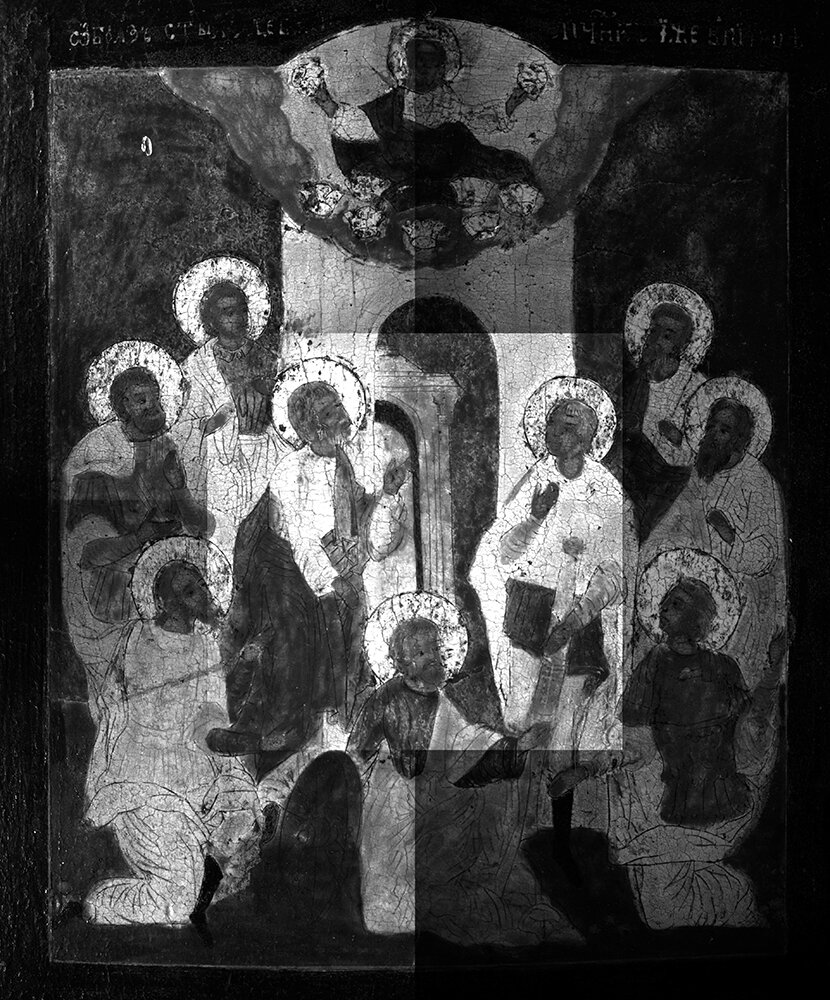}\,
\includegraphics[width=0.13\textwidth]{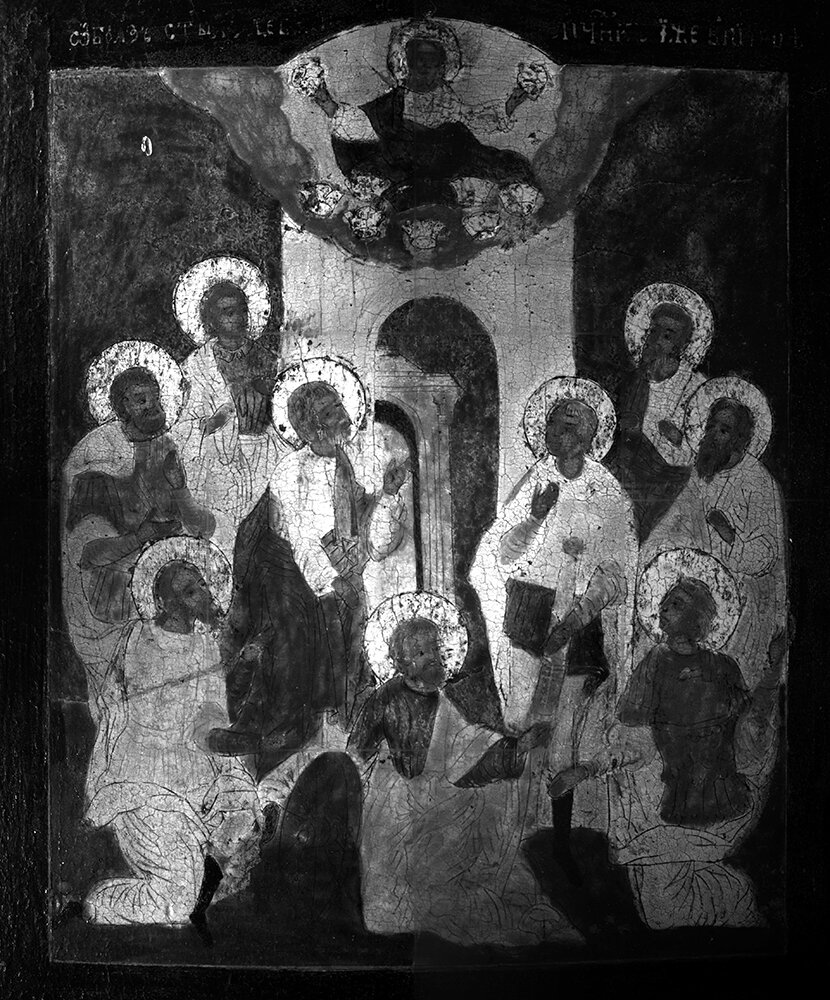}\,
\includegraphics[width=0.10\textwidth,trim=0 0 10cm 0,clip=true]{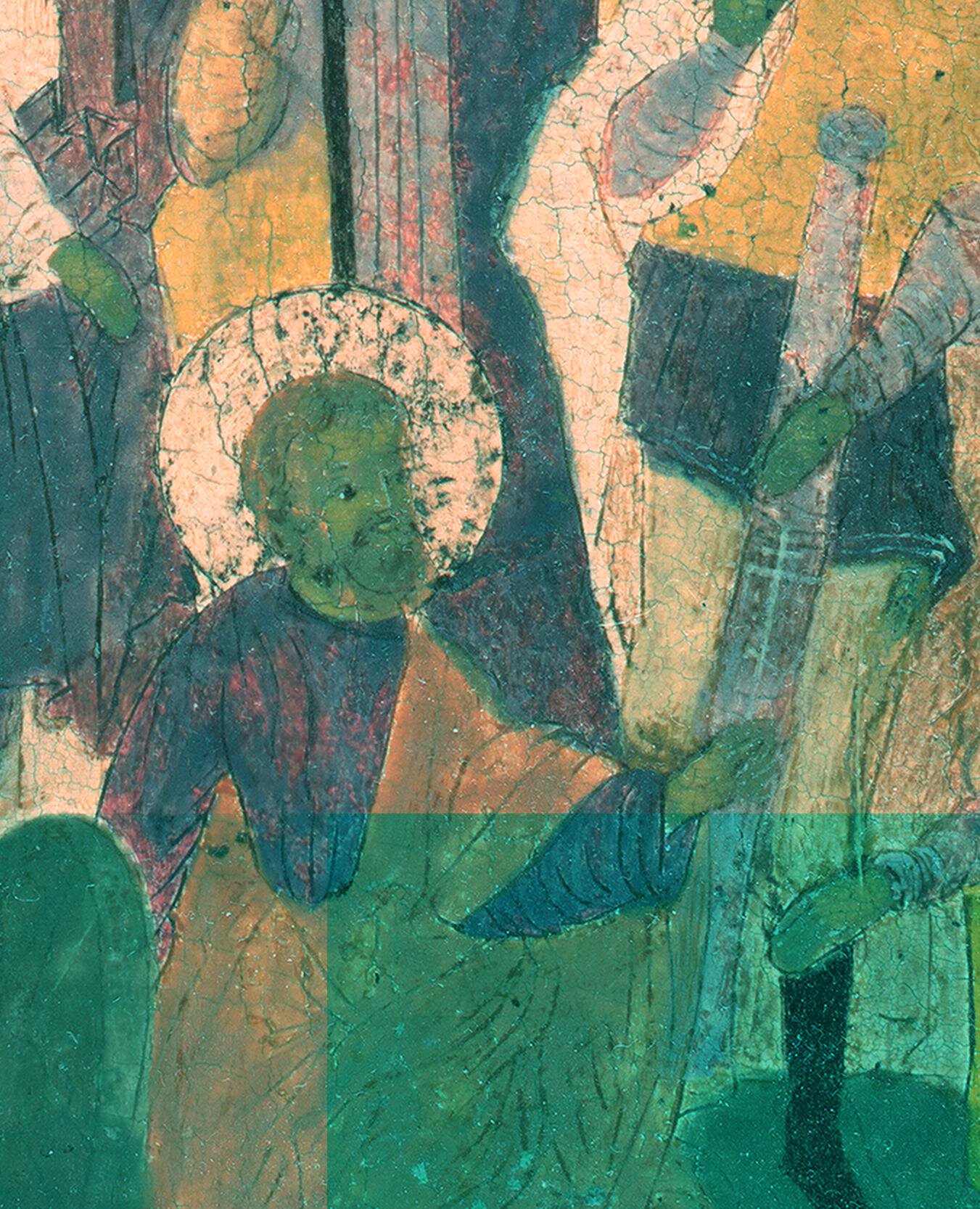}\,
\includegraphics[width=0.10\textwidth,trim=0 0 10cm 0,clip=true]{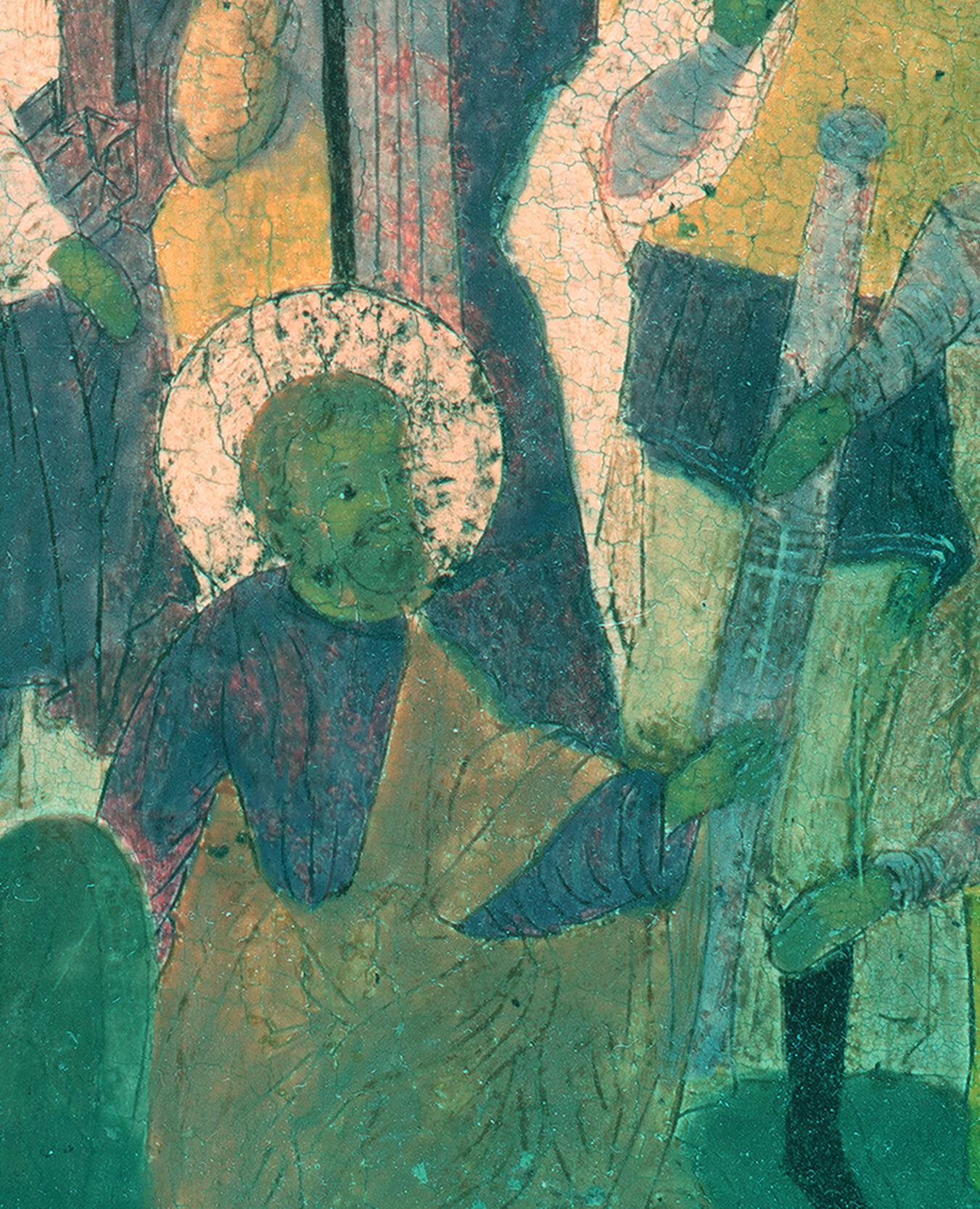}
\caption{Light balancing in IR falsecolor imaging. IR before and after osmosis (left and middle-left), IR-R-G before and after osmosis (middle-right and right, detail of 18 MP, 721 s.\ \cite{Parisotto2018})}
\label{fig: ir falsecolor}
\end{figure}

In CH imaging there exists other situations in which multiple images are sampled under different light conditions, which may introduce errors and affect the quality of the final image to inspect.
Again, the use of an efficient numerical scheme of the osmosis model \eqref{eq: osmosis} makes the task of homogenizing the light differences possible in a reasonable time.
In Figs.\ \ref{fig: uvf} and \ref{fig: ir falsecolor} we report a case study for the UV fluorescence inspection and the falsecolor imaging (IR-Red-Green) of a Russian icon.

\begin{figure}[tbhp]
\centering
\includegraphics[width=0.13\textwidth,trim=0 2cm 0 0,clip=true]{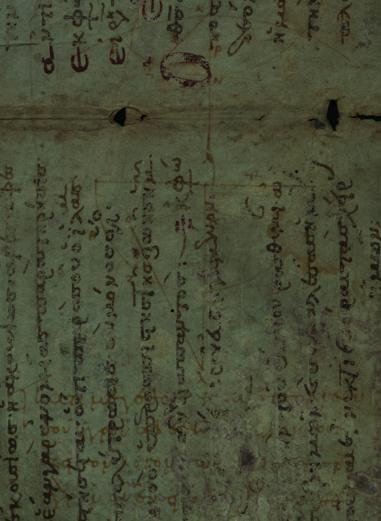}\,
\includegraphics[width=0.13\textwidth,trim=0 2cm 0 0,clip=true]{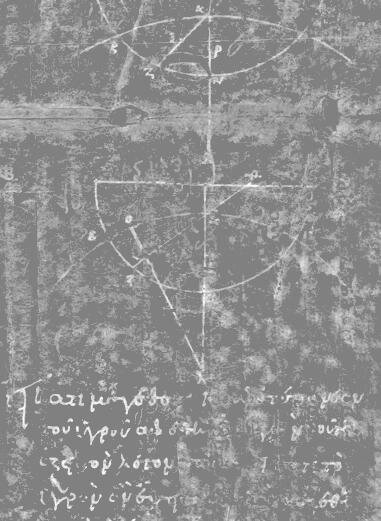}\,
\includegraphics[width=0.13\textwidth,trim=0 2cm 0 0,clip=true]{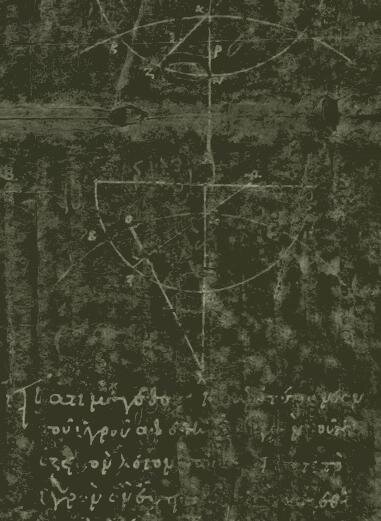}
\caption{Data integration: \emph{Archimedes Palimpsest} (detail, 2 MP, 137 s. \cite{Parisotto2018}). Visible parchment (left), hidden text (center), data fusion (right).}
\label{fig: archi}
\end{figure}

In Fig.\ \ref{fig: archi} we report a detail of the \emph{Archimedes Palimpsest}, a X century copy of the works by Archimedes overwritten in XIII century. 
Here we used the property of \eqref{eq: osmosis} to converge to a rescaled version of $v$, encoded in $\bf{d}$, to fuse the original written text with the current colours.

\section{Conclusion}
\label{sec: conclusions}
In this work, we highlighted the dual-mode MWIR analysis for the restorations of frescoes \cite{Daffara2018}. 
Also, we tested an efficient numerical approach of the osmosis equation \cite{Vogel2013,Weickert2013} for a variety of CH multi-spectral applications, from balancing lights in multi-spectral images to data fusion on Archimedes Palimpsest \cite{Parisotto2018}.

\paragraph*{Data statement.} The \emph{Monocromo} data are sensitive: access subjected to \emph{Soprintendenza Castello Sforzesco}, Milan. The Archimedes Palimpsest\footnote{http://openn.library.upenn.edu/Data/0014/ArchimedesPalimpsest} is released under CC-BY 3.0 license.


\end{document}